\documentclass[12pt]{article}
\usepackage[final]{epsfig}
\usepackage{graphics}
\usepackage{amsmath}
\usepackage{amsfonts}
\usepackage{latexsym}
\usepackage{amssymb}
\usepackage{graphicx}
\usepackage{epstopdf}

\newtheorem{lemma}{Lemma}[section]
\newtheorem{lemma-def}[lemma]{Lemma-Definition}
\newtheorem{proposition}[lemma]{Proposition}
\newtheorem{remark}[lemma]{Remark}

\newtheorem{theorem}[lemma]{Theorem}

\newtheorem{conjecture}[lemma]{Conjecture}
\newtheorem{problem}[lemma]{Problem}

\begin{document}
\newcommand{\eps}{{\varepsilon}}
\newcommand{\proofend}{$\Box$\bigskip}
\newcommand{\C}{{\mathbf C}}
\newcommand{\Q}{{\mathbf Q}}
\newcommand{\R}{{\mathbf R}}
\newcommand{\Z}{{\mathbf Z}}
\newcommand{\RP}{{\mathbf {RP}}}

\def\proof{\paragraph{Proof.}}


\newcommand{\marginnote}[1]
{
}

\newcounter{bk}
\newcommand{\bk}[1]
{\stepcounter{bk}$^{\bf BK\thebk}$%
\footnotetext{\hspace{-3.7mm}$^{\blacksquare\!\blacksquare}$
{\bf BK\thebk:~}#1}}

\newcounter{st}
\newcommand{\st}[1]
{\stepcounter{st}$^{\bf ST\thest}$%
\footnotetext{\hspace{-3.7mm}$^{\blacksquare\!\blacksquare}$
{\bf ST\thest:~}#1}}

\newcounter{dg}
\newcommand{\dg}[1]
{\stepcounter{dg}$^{\bf DG\thebk}$%
\footnotetext{\hspace{-3.7mm}$^{\blacksquare\!\blacksquare}$
{\bf DG\thebk:~}#1}}


\title {Geodesics on an ellipsoid in Minkowski space}
\author
{Daniel Genin,\thanks{
Department of Mathematics,
Pennsylvania State University, University Park, PA 16802, USA;
e-mail: \tt{genin@math.psu.edu}
}
\, Boris Khesin\thanks{
Department of Mathematics,
University of Toronto, Toronto, ON M5S 2E4, Canada;
e-mail: \tt{khesin@math.toronto.edu}
}
\, and Serge Tabachnikov\thanks{
Department of Mathematics,
Pennsylvania State University, University Park, PA 16802, USA;
e-mail: \tt{tabachni@math.psu.edu}
}
\\
}
\date{May 1, 2007}
\maketitle
\begin{abstract}
We describe the geometry of geodesics on a Lorentz ellipsoid: give explicit formulas for
the first integrals (pseudo-confocal coordinates), curvature, geodesically equivalent Riemannian metric, the invariant area-forms on the time- and space-like geodesics and invariant 1-form on the space of null geodesics. We prove a Poncelet-type theorem for null geodesics on the ellipsoid: if such a geodesic close up after several oscillations in the ``pseudo-Riemannian belt", so do all other  null geodesics on this ellipsoid. 
\end{abstract}

\section{Introduction} \label{intro}

The geodesic flow on the ellipsoid in Euclidean space is a classical example of a completely integrable dynamical system whose study goes back to Jacobi and Chasles. We refer to \cite{A-G,Kn,Mo1} 
 for a modern treatment of this subject; see also \cite{AF,Au,CS,DR,Pr,Ve} for a sampler of recent work. In a recent paper \cite{Kh-T}, we considered ellipsoids in pseudo-Euclidean spaces of arbitrary signatures and extended, with appropriate adjustments, the theorem on complete integrability of the geodesic flow to this case. We also defined pseudo-Euclidean billiards and proved that the billiard map inside an ellipsoid in pseudo-Euclidean space is completely integrable.

This note is devoted to the ``case study" of geodesics on an ellipsoid in three dimensional space with the metric $dx^2+dy^2-dz^2$.  Recall that the pseudo-sphere $x^2+y^2-z^2=-1$ in Minkowski space provides a famous example of a Riemannian metric of constant negative curvature, a model of the hyperbolic plane, and this gives a further motivation to investigate quadratic surfaces in Minkowski space.

We shall study an ellipsoid
\begin{equation} \label{triell}
\frac{x^2}{a}+\frac{y^2}{b}+\frac{z^2}{c}=1,\quad a,b,c>0,
\end{equation}
and we assume that the general position condition $a> b$ holds.
The induced metric $\langle\,,\rangle$ on the ellipsoid degenerates along the two curves
\begin{equation} \label{tropic}
z=\pm c\sqrt{\frac{x^2}{a^2} + \frac{y^2}{b^2}}
\end{equation}
that will be referred to as the ``tropics". The induced metric \ is Riemannian in the ``polar caps" and Lorentz in the ``equatorial belt" bounded by the tropics. 
We shall see in Section \ref{equi} that the Gauss curvature of the ellipsoid is everywhere negative (it equals $-\infty$ on the tropics).
At every point of the equatorial belt, one has two null directions of the Lorentz metric; on the tropics, these directions merge together.

Every geodesic curve $\gamma(t)$ on the ellipsoid (\ref{triell}) is of one of the three types: space-like (positive energy $\langle\gamma',\gamma'\rangle>0$), time-like (negative energy $\langle\gamma',\gamma'\rangle<0$) or light-like (zero energy or null $\langle\gamma',\gamma'\rangle=0$). The last two types exist only in the equatorial belt. In fact, a whole new phenomenon which we are dealing with in the Lorentz case,
as compared to the Euclidean one, is the presence of null geodesics, which
``separate" the space-like and time-like ones. This makes the
pseudo-Riemannian geometry of the problem quite peculiar, and rather
different from its Riemannian counterpart.

Let us summarize some relevant results from \cite{Kh-T}, specialized to the three dimensional case. Include the ellipsoid into a {\it pseudo-confocal} family of quadrics $M_{\lambda}$ given by the equation
\begin{equation} \label{pconf}
\frac{x^2}{a+\lambda}+\frac{y^2}{b+\lambda}+\frac{z^2}{c-\lambda}=1.
\end{equation}
It is shown in \cite{Kh-T} that: 
\begin{enumerate}
\item Through every generic point $Q(x,y,z)$ in space, there pass either three or one pseudo-confocal quadrics. In the former case, two of the quadrics have the same topological type, and the quadrics are pairwise orthogonal at point $Q$. If $Q$ is a generic point of the ellipsoid (\ref{triell}) then there exists another pseudo-confocal ellipsoid and a pseudo-confocal hyperboloid of one sheet, passing through $x$.
 
\item A generic space- or time-like line $\ell$  is tangent to either  two or no pseudo-confocal quadrics. In the former case, the tangent planes to these quadrics at the tangency points with $\ell$ are pairwise orthogonal. In particular, a generic line, tangent to  the ellipsoid (\ref{triell}), is tangent to another  pseudo-confocal quadric.  

\item The tangent lines to a given space-like or time-like geodesic on the ellipsoid (\ref{triell}) forever remain tangent to a fixed pseudo-confocal quadric; the respective value of $\lambda$ in (\ref{pconf}) can be considered as an integral of the geodesic flow. 
\end{enumerate}
 
The spaces of space- or time-like lines in pseudo-Euclidean space, and more generally, the spaces of space- or time-like geodesics in a pseudo-Riemannian manifold, carry symplectic structures. These structures are obtained by symplectic reduction: restrict the canonical symplectic structure of the (co)tangent bundle to the unit energy hypersurface $\langle v,v\rangle =\pm 1$ and quotient out the one-dimensional kernel of this restriction. In the case of a pseudo-Riemannian (Lorentz) surface, we obtain an area form on the set of space- and time-like geodesics.\footnote{It is shown in \cite{Kh-T} that the space of light-like geodesics has a contact structure; it is trivial in the case of a surface, as  light-like geodesics form a one-dimensional set, and we shall not use this structure in the present paper.}  
 
Now, to billiards.  In general, the billiard dynamical system in a pseudo-Riemannian manifold with a smooth boundary describes the motion of a free mass-point (``billiard ball"). The point moves along a geodesic with constant energy until it hits the boundary where the elastic reflection occurs: the normal component of the velocity instantaneously changes sign whereas the tangential component remains the same, see \cite{Ta5,Ta2} for general information. The billiard reflection is not defined at the   points where the normal vector is tangent to the boundary. The {\it billiard ball map} acts on oriented geodesics and takes the incoming trajectory of the  billiard ball to the outgoing one. This map preserves the type of a geodesic (space-, time-, or light-like). Furthermore, the billiard ball map, acting on the space- or time-like geodesics, preserves the symplectic structure (area form, in two-dimensional case) described above.\footnote{And, acting on light-like geodesics, preserves the contact structure.}

It is proved in \cite{Kh-T} that the billiard ball map with respect to a quadric in Minkowski space is integrable as well: if an incoming space- or time-like billiard trajectory is tangent to two pseudo-confocal quadrics then the outgoing trajectory remains tangent to the same two quadrics.
 
\begin{figure}[hbtp]
\centering
\includegraphics[width=2in]{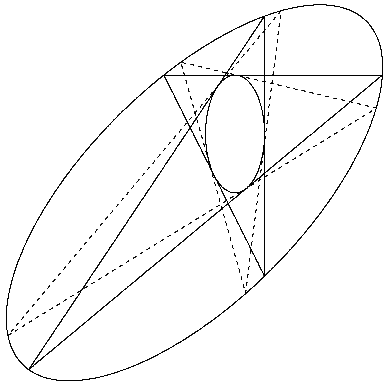}
\caption{Poncelet Porism}
\label{Ponc}
\end{figure}

In general, complete integrability of a billiard system implies   configuration theorems for closed billiard trajectories, see, e.g., \cite{Ta5,Ta2}. In particular, the integrability of the billiard inside an ellipse implies the classical Poncelet Porism depicted in figure \ref{Ponc}: let $\gamma\subset \Gamma$ be two nested ellipses and let a point of $\Gamma$ be a vertex of an $n$-gon inscribed in $\Gamma$ and circumscribed about $\gamma$; then every point of $\Gamma$ is a vertex of such an $n$-gon, see \cite{B-K-O-R}.

In Section \ref{Poncelet} we prove a Poncelet-style closure theorem involving null geodesics on the ellipsoid: if some chain of such geodesics in the equatorial belt  closes up after $n$  oscillations between the tropics then so does every such chain.

It is proved in \cite{Kh-T} that the billiard inside an ellipse in the Lorentz plane is completely integrable as well. Restricted to light-like lines, the billiard ball map in an  oval becomes a circle map depending only on the two null directions. In Section \ref{rig} we discuss the problem when this map is conjugated to a rotation.
 
\section{Joachimsthal  integral} \label{Jint}

The geodesic flow on the triaxial ellipsoid in Euclidean space possesses the classical Joachimsthal  integral. In this section we describe its pseudo-Euclidean version. This integral says analytically   that the lines tangent to a geodesic curve on the ellipsoid are tangent to a fixed pseudo-confocal quadric: in this dimension, all first integrals of the geodesic flow, restricted to a constant energy hypersurface,  are functionally dependent. 

Choose a Minkowski normal at point $(x,y,z)$ of the ellipsoid  (\ref{triell}) as follows: 
$$
N(x,y,z)=\left( \frac{x}{a},\frac{y}{b}, -\frac{z}{c} \right).
$$
Then 
$$
\langle N(x,y,z),N(x,y,z)\rangle=\frac{x^2}{a^2}+\frac{y^2}{b^2} -\frac{z^2}{c^2}.
$$
We see that the normal is space-like in the equatorial belt, time-like in the polar caps and null on the tropics. 

Let $(u,v,w)$ denote a tangent vector to the ellipsoid.

\begin{proposition} \label{Joachimsthal}
The following function is an integral of the geodesic flow:
\begin{equation} \label{Joa}
J=\left( \frac{x^2}{a^2}+\frac{y^2}{b^2}-\frac{z^2}{c^2}\right) \left( \frac{u^2}{a}+\frac{v^2}{b}+\frac{w^2}{c}\right).
\end{equation}
\end{proposition}

\proof
Parameterized geodesics $\gamma(t)$ are extrema of the functional 
$$
\int \langle\dot \gamma(t), \dot \gamma(t)\rangle\ dt;
$$
they satisfy the Euler-Lagrange equation
$
\ddot \gamma(t)=\lambda(t) N(\gamma(t))
$
where  $\lambda(t)$ is a Lagrange multiplier.

Let $\gamma(t)$ be a geodesic with $\dot\gamma(t) =(u,v,w)$. Then $\ddot \gamma(t)=\lambda(t) N(\gamma(t))$, and hence 
$$
\lambda(t)=\frac{\langle \ddot \gamma(t),N(\gamma(t))\rangle}{\langle N(\gamma(t)),N(\gamma(t))\rangle}.
$$
Since $\langle \dot \gamma(t), N(\gamma(t))\rangle=0$, we have: 
$$
\langle\ddot \gamma(t), N(\gamma(t))\rangle=- \langle\dot \gamma(t), \dot N(\gamma(t))\rangle= 
- \langle(u,v,w), \left( \frac{u}{a},\frac{v}{b}, -\frac{w}{c}\right)\rangle=-\left( \frac{u^2}{a}+\frac{v^2}{b}+\frac{w^2}{c}\right),
$$
and it follows that 
\begin{equation} \label{lambda}
\lambda(t)=-\frac{\frac{u^2}{a}+\frac{v^2}{b}+\frac{w^2}{c}}{\frac{x^2}{a^2}+\frac{y^2}{b^2}-\frac{z^2}{c^2}}
\end{equation}
(the denominator vanishes precisely on the tropics).

We want to prove that $\dot J =0$. Indeed,
$$
2\dot J = \left( \frac{xu}{a^2}+\frac{yv}{b^2}-\frac{zw}{c^2}\right) \left( \frac{u^2}{a}+\frac{v^2}{b}+\frac{w^2}{c}\right) + \left( \frac{x^2}{a^2}+\frac{y^2}{b^2}-\frac{z^2}{c^2}\right) \left( \frac{u\dot u}{a}+\frac{v\dot v}{b}+\frac{w\dot w}{c}\right)
$$
$$=
\left( \frac{xu}{a^2}+\frac{yv}{b^2}-\frac{zw}{c^2}\right) \left( \frac{u^2}{a}+\frac{v^2}{b}+\frac{w^2}{c}\right) + \lambda \left( \frac{x^2}{a^2}+\frac{y^2}{b^2}-\frac{z^2}{c^2}\right)
 \left( \frac{xu}{a^2}+\frac{yv}{b^2}-\frac{zw}{c^2}\right)=0, 
$$
the last equality due to (\ref{lambda}).
\proofend


\section{Describing the geodesics} \label{descgeo}

\subsection{Global behavior of geodesics} \label{glob}

As we mentioned earlier, the metric in the polar caps is Riemannian, while  
in the equatorial belt it is Lorentzian. A geodesic in a polar cap crosses it  from tropic to tropic.

Consider the equatorial belt. The light-like geodesics traverse the belt from one tropic to another. There are two kinds of them depending on the slope, positive or negative; we call these null geodesics {\it right} and {\it left}, respectively. The time-like geodesics are squeezed between like-like ones and hence also go from tropic to tropic. 

\begin{figure}[hbtp]
\centering
\includegraphics[width=4.5in]{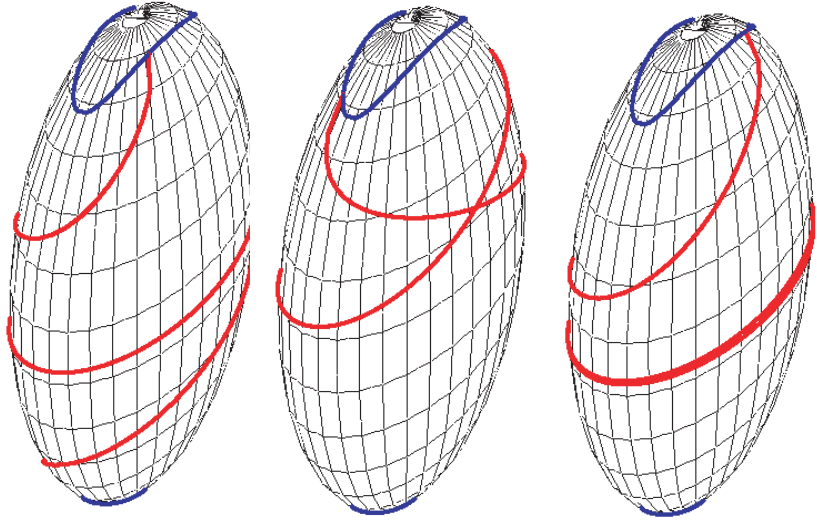}
\caption{Space-like geodesics in the equatorial belt}
\label{eqspace}
\end{figure}

Let us discuss the behavior of  space-like geodesics. 
The equator is a closed geodesic. Figure \ref{eqspace} depicts two other space-like geodesics (first two pictures); the first one intersects the equator  and proceeds to the other tropic, and the other geodesic starts at a tropic and returns to the same tropic without reaching the equator. These two geodesics represent generic behavior as the next lemma shows.

\begin{lemma} \label{turn}
Let $\gamma(t)=(x(t),y(t),z(t))$ be a geodesic in the Northern part of the equatorial belt, i.e., $z(t)>0$. Then $\ddot z(t)>0$. In the Southern part of the equatorial belt, $\ddot z(t)<0$
\end{lemma}

\proof
Since $\ddot \gamma(t)=\lambda(t) N(\gamma(t))$, we have: $\ddot z(t)=-\lambda(t) z(t)/c$.  According to (\ref{lambda}), $\lambda<0$ in the equatorial belt, hence $\ddot z(t)>0$
for $z>0$ (i.e.,  in the Northern part) and $\ddot z(t)<0$
for $z<0$ (in the Southern part). 
\proofend

It follows that a geodesic in the Northern part of the equatorial belt is convex downwards, and in the Southern part of the equatorial belt is convex upwards. 
It also follows from the proof that the equator, as a closed geodesic, is exponentially unstable. 

An intermediate position between the two types of generic space-like geodesics is occupied by the geodesics that start at, say, Northern tropic and monotonically descend to the equator having the latter as a limit cycle, see figure \ref{eqspace}, on the right.

\subsection{Local behavior of geodesics near the tropics} \label{near}

The behavior of geodesics near tropics is described by the following

\begin{proposition} \label{tropdir}
A geodesic on the ellipsoid  (\ref{triell}) may reach a tropic only in the null direction.\end{proposition}

\proof
For light-like geodesics, there is nothing to prove, and a time-like geodesic is confined in the wedge between two null directions which merge together on a tropic. It remains to consider unit energy space-like geodesics. 

According to Proposition \ref{Joachimsthal}, one has, along a geodesic:
$$
\left( \frac{x^2}{a^2}+\frac{y^2}{b^2}-\frac{z^2}{c^2}\right) \left( \frac{u^2}{a}+\frac{v^2}{b}+\frac{w^2}{c}\right)=const,\ \ u^2+v^2-w^2=1.
$$
As the point approaches a tropic, the first factor on the left hand side of the first equality goes to zero, and hence the second factor blows up. The direction of the geodesic does not change if we rescale the tangent vector:
$$
(\bar u,\bar v,\bar w)=\mu (u,v,w),\ \ \mu=\sqrt{\frac{x^2}{a^2}+\frac{y^2}{b^2}-\frac{z^2}{c^2}}
$$
maintaining ${\bar u}^2+{\bar v}^2+{\bar w}^2 =O(1)$. As $(x,y,z)$ approaches the tropic, we have 
$$
{\bar u}^2+{\bar v}^2-{\bar w}^2 = \mu^2 \to 0,
$$
therefore the geodesic has the null direction.
\proofend

It follows from Proposition \ref{tropdir} that the null lines tangent to the ellipsoid along the tropics are tangent to infinitely many pseudo-confocal quadrics. Indeed, infinitely many geodesic hit the tropics at any given point in the null direction, and the tangent lines to each geodesic are tangent to a fixed pseudo-confocal quadric, see Section \ref{intro}. In contrast, we have the following proposition which provides an alternative proof of Proposition \ref{tropdir}.

\bigskip
\noindent{\bf Proposition \ref{tropdir}$'$} \label{none}
{\it 
Let $P$ be a point of a tropic and $\ell$ a space-like line tangent to the 
ellipsoid at $P$. Then $\ell$ is not tangent to any pseudo-confocal 
quadric.
}

\proof
Assume that $\ell$ is tangent to $M_{\lambda}$ at point $Q$. Denote by $N$ a normal vector to $M_0$ at $P$ and by $\eta$ a normal vector to $M_{\lambda}$ at $Q$. Then $N$ lies in the tangent plane $T_P M_0$. The restriction of the ambient metric to this plane is degenerate, and $N$ spans the kernel of this restriction.

According to the general theory, see Section \ref{intro}, $\eta$ is orthogonal to $N$, hence $\eta$ lies in $N^{\perp}=T_P M_0$. Since $\eta$ is also orthogonal to $\ell$, the normals $\eta$ and $N$ are collinear. Therefore the tangent planes $T_P M_0$ and $T_Q M_{\lambda}$ coincide. 

\begin{figure}[hbtp]
\centering
\includegraphics[width=2in]{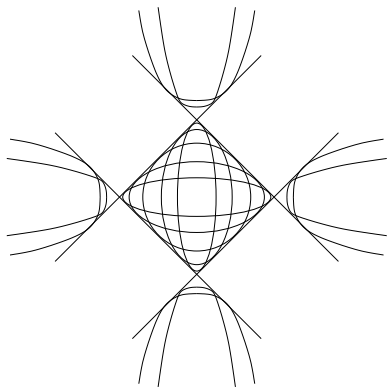}
\caption{Pseudo-confocal family of conics $\frac{x^2}{1+\lambda}+\frac{y^2}{1-\lambda}=1$}
\label{family}
\end{figure}

The projection of pseudo-confocal family along the line $\ell$ is a pseudo-confocal family of conics in the Lorentz plane, see \cite{Kh-T}. It follows that the projections of $M_0$ and $M_{\lambda}$ are tangent to each other. However no two conics in a pseudo-confocal family are tangent (see figure \ref{family}),  a contradiction. 
\proofend

Proposition \ref{tropdir}  is a manifestation of a more general phenomenon which holds for  Lorentz surfaces. Suppose, on an analytic Lorentz surface $S$, the Lorentz metric changes into Riemannian (i.e., the two null direction coincide) along a smooth curve $\Gamma$. We call a point of $\Gamma$ {\it regular} if the null direction at this point is transversal to the curve $\Gamma$ itself.

\begin{figure}[hbtp]
\centering
\includegraphics[width=1.5in]{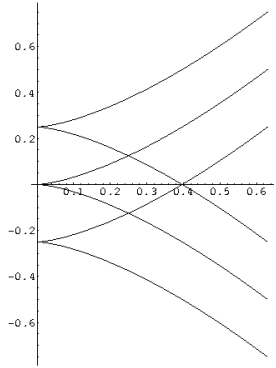}
\caption{Family of cusps}
\label{cusps}
\end{figure}

\begin{proposition} \label{normalform}
In a neighborhood of a regular point the Lorentz metric on the surface $S$ is conformally equivalent to the metric $dy^2-x\,dx^2$ in a neighborhood of the origin.
The null geodesics  on $S$ near this point are diffeomorphic to the family of cusps
$y=x^{3/2}+C$, see figure \ref{cusps}.
\end{proposition}

\proof

Let $a(x,y)dx^2+ b(x,y)dxdy+c(x,y)dy^2$ be a Lorentz metric on $S$. Then the equation of null geodesics on $S$ is $a(x,y)+ b(x,y)y'+c(x,y)(y')^2=0$, an implicit differential equation of the first order. The condition of transversality at a regular point of $\Gamma$
means that we consider a so-called ``regular singular point" of the implicit differential equation $F(x,y,y')=0$.
According to the theorem of Cibrario (see a discussion in \cite{Ar}) the normal form of this differential equation at a regular singular point is $(y')^2=x$, which is the equation of null geodesics in the metric $dy^2-x\,dx^2$. Thus the null geodesics is the family of cusps  described above.

The family of null geodesics fixes the conformal class of the Lorentz metric. Due to analyticity, this metric uniquely extends beyond the curve $\Gamma$, into the Riemannian domain of the neighborhood of the regular point.
\proofend

\begin{remark}
{\rm
Note that every geodesic for the metric $dy^2-x\,dx^2$ which hits the ``tropic" $\Gamma=\{x=0\}$, does it with a horizontal velocity. Indeed, consider a non-vertical space-like geodesic with the velocity vector $(u,v)$ of the unit Lorentz length: $v^2-x\,u^2=1$. Note that $v^2\not=1$ for any $x\not=0$.
Since the metric is invariant with respect to $y$-translations, the velocity vector along any geodesic conserves its $v$-component (in addition to the
conservation of its Lorentz length). Then, at the moment  of ``impact" with the tropic, the $u$-component has to become infinite: $u^2=(v^2-1)/x\to\infty$ as $x\to 0$, that is, the velocity becomes horizontal on the tropic.

Although, for a metric which is only conformally equivalent to this normal form, one cannot use the invariance of the $v$-component, one can bound it above and below, and hence the $u$-component has to become infinite anyway, i.e., the same conclusion on the horizontality of the velocity holds, due to its robustness.
This consideration implies Proposition \ref{tropdir} for the ellipsoid, once one checks that the null directions are everywhere transversal to the tropics, which is straightforward. 
}
\end{remark}

\subsection{Intersections  with pseudo-confocal quadrics}\label{inter}

\begin{figure}[hbtp]
\centering
\includegraphics[width=4in]{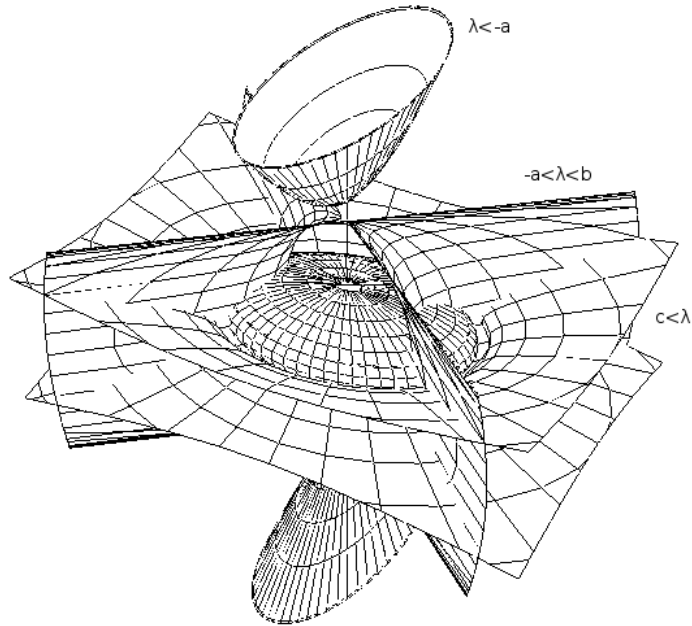}
\caption{Four pseudo-confocal quadrics, along with the initial ellipsoid}
\label{Confquad}
\end{figure}

\begin{figure}[hbtp]
\centering
\includegraphics[width=2.5in]{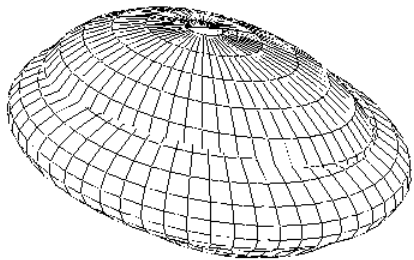}
\caption{A pseudo-confocal ellipsoid intersecting the initial one}
\label{Confell}
\end{figure}

The topology of a pseudo-confocal quadric (\ref{pconf}) depends on the position of $\lambda$ relative the three numbers $-a, -b$ and $c$: if $\lambda< -a$ then $M_{\lambda}$ is a hyperboloid of two sheets, if $-a <\lambda < -b$ then $M_{\lambda}$ is a hyperboloid of one sheet, if $-b <\lambda < c$ then $M_{\lambda}$ is an ellipsoid, and if $c <\lambda$ then $M_{\lambda}$ is again a hyperboloid of one sheet. Only the second and the third kinds intersect the original ellipsoid, $M_0$. See figure \ref{Confquad}
for all five quadrics and figure \ref{Confell} for an ellipsoid $M_{\lambda}$ with $-b <\lambda < c$ intersecting $M_0$.

The intersection curve $\Gamma_{\lambda}=M_0 \cap M_{\lambda}$ is given by the equation
\begin{equation} \label{intcurv}
\frac{x^2}{a(a+\lambda)}+\frac{y^2}{b(b+\lambda)}=\frac{z^2}{c(c-\lambda)}.
\end{equation}
In the limit  $\lambda \to 0$, this becomes the equation of the tropic (\ref{tropic}). As $\lambda \to c$, the curve $\Gamma_{\lambda}$ tends to the equator. 
The projection of the intersection curve $\Gamma_{\lambda}$ on the $(x,y)$-plane is a conic
$$
\frac{a+c}{a(a+\lambda)} x^2 + \frac{b+c}{b(b+\lambda)} y^2 =1.
$$
If $-a <\lambda < -b$ then this conic is a hyperbola, and if $-b <\lambda < c$ it is an ellipse. 

Fix a value of $\lambda$. Given a generic point $P$ of the ellipsoid $M_0$, the number of tangent lines from $P$ to the conic  $T_P M_0 \cap M_{\lambda}$ may equal 2  or 0, and the curve $\Gamma_{\lambda}$ separates these two domains, say, $U_{\lambda}$ and $V_{\lambda}$ (the number of tangents equals 1 for $P \in \Gamma_{\lambda}$). 

Consider a geodesic $\gamma$ on $M_0$. The lines tangent to $\gamma$ are tangent to some pseudo-confocal quadric $M_{\lambda}$. Therefore $\gamma$ is confined to $U_{\lambda}$ and, at any point $P \in U_{\lambda}$, the geodesic $\gamma$ may have 
only two directions:  the tangent  directions from $P$ to the conic  $T_P M_0 \cap M_{\lambda}$. The geodesic $\gamma$ cannot intersect the boundary $\Gamma_{\lambda}$ but can touch it. 

Note that the equator plays a special role: no line tangent to the equator is tangent to any pseudo-confocal quadric (except $M_0$). Indeed, such a line lies in the horizontal plane, and the trace of the pseudo-confocal family (\ref{pconf}) in the horizontal plane is the confocal family of conics. But confocal conics have no common tangents. 

\subsection{Reflection of  geodesics off the tropics}\label{reflspace}

Define ``reflection" of a geodesic $\gamma$ off a tropic as the geodesic  tangent to the same pseudo-confocal quadric as $\gamma$ or, equivalently, having the same value of the Joachimsthal  integral (normalizing the energy to $\pm 1$). According to Proposition \ref{tropdir},  a geodesic and a reflected one have the null direction at the impact point and hence are tangent to each other.

\begin{figure}[hbtp]
\centering
\includegraphics[width=3in]{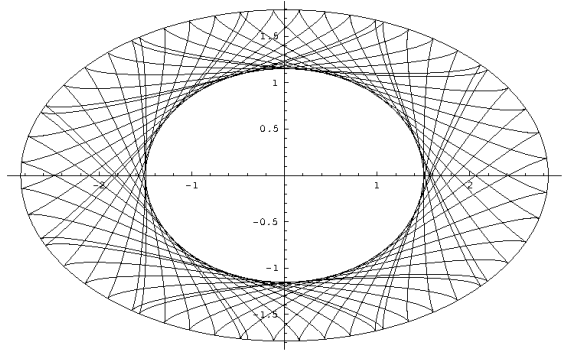}
\caption{Reflected geodesics in a polar cap}
\label{polarbil}
\end{figure}

\begin{figure}[hbtp]
\centering
\includegraphics[width=3in]{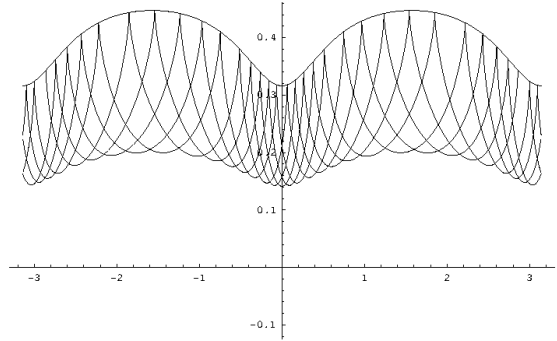}
\caption{Reflected geodesics in the equatorial belt (unfolded)}
\label{eqbil}
\end{figure}

\begin{proposition} \label{envel}
1). Let $\gamma$ be a geodesic in a polar cap tangent to a curve $\Gamma_{\lambda}=M_0 \cap M_{\lambda}$, and let $\gamma_1$ be the reflection of $\gamma$ in the tropic. Then $\gamma_1$ is also tangent to $\Gamma_{\lambda}$ (see figure \ref{polarbil}). \hfill \break
2). Let  $\gamma$ be a space-like geodesic in the equatorial belt that does not intersect the equator,  and let $\gamma_1$ be the reflection of $\gamma$ in the tropic. Then $\gamma$ is tangent to a  curve $\Gamma_{\lambda}$, separating $\gamma$ from the equator, and $\gamma_1$ is also tangent to $\Gamma_{\lambda}$ and therefore disjoint from the equator, see figure \ref{eqbil}.
\end{proposition}

\proof
 Since the geodesic $\gamma$ is tangent to $\Gamma_{\lambda}$, the tangent  lines to $\gamma$ are tangent to $M_{\lambda}$.  By definition, the tangent  lines to $\gamma_1$ are tangent to $M_{\lambda}$ as well, and hence $\gamma_1$ will touch $\Gamma_{\lambda}$, and the first claim follows.

Likewise, the tangent lines to  a space-like geodesic $\gamma$ in the equatorial belt are tangent to a pseudo-confocal geodesic $M_{\lambda}$.   Since $\gamma$ does not reach the equator, the curve $\Gamma_{\lambda}$ separates the tropic and the equator. Arguing as above, the reflected geodesic $\gamma_1$ is also tangent to $\Gamma_{\lambda}$ and therefore disjoint from the equator.
\proofend

\subsection{Geodesics on degenerate ellipsoids}\label{pancakes}

It is useful   to visualise the degenerations of the ellipsoid (\ref{triell}) 
into a two-sided flat surface, ``a pancake." First consider the limit  
$c\to 0$, under which the ellipsoid gets squeezed to an ellipse  with semi-axes $a$ and $b$ in the Euclidean plane $\{z=0\}$. Then the (Lorentz) equatorial belt disappears, 
whereas each (Riemannian) polar cap becomes the interior of an ellipse. 
The geodesics in the polar caps
become straight lines in the limit, and their reflection off the tropics becomes
the billiard reflection in the ellipse; the Joachimsthal integral describes the confocal
ellipse to which a billiard trajectory remains tangent.

Compare this with the other degeneration, $b\to 0$. In the latter case
the ellipsoid becomes (the interior of) an ellipse in the Lorentz plane $\{y=0\}$. The polar caps 
(and tropics) become squeezed to the  arcs of the ellipse $z>0$ and $z<0$ in this plane. The equatorial belt becomes the double cover of the  interior of the ellipse. The geodesics also become  straight lines  in the Lorentz plane. The null geodesics have two prescribed slopes in this plane: $x=\pm z$. We will see that the dynamics of the null geodesics in the ellipsoid is an interesting extension of the corresponding Lorentz billiard inside the ellipse, restricted to oriented null lines.


\section{Area form on the space of geodesics} \label{area}

In this section we compute the area form $\omega$ on the space of time-like geodesics, described in Section \ref{intro}. We also define a 1-form on the space of light-like geodesics which will play the central role in the next section.

Let us characterize a time-like  or a light-like   geodesic by its intersection with the equator of the ellipsoid (\ref{triell}). The equator is parameterized as $Q(t)=(\sqrt{a} \cos t, \sqrt{b} \sin t, 0)$ for $t\in\R/2\pi\Z$. 
Let a geodesic make (Euclidean) angle $\alpha$ with the equator. Then $(t,\alpha)$ are coordinates in the space of geodesics  intersecting the equator.  The null directions correspond to $\alpha=\pi/4$ and $3\pi/4$, and for time-like geodesics, $\pi/4 < \alpha< 3\pi/4$. 
We set 
$$
f(t)=\sqrt{a\sin^2t+b\cos^2t},\ \  \tau=\frac{1}{\sqrt{\tan^2\alpha -1}}.
$$
Note that for $\alpha=\pi/2$, the value of $\tau$ is well-defined: $\tau=0$.
The following proposition describes the Joachimsthal  integral $J$ and the area form $\omega$  in terms of the $(t,\alpha)$ coordinates.

\begin{proposition} \label{onneg}
One has:
\begin{equation} \label{JandO}
 J(t,\alpha)= \frac{c\tau^2 + f^2(t)(1+\tau^2)}{abc},\ \ \omega=f(t) d\tau \wedge dt.
\end{equation}
\end{proposition}

\proof For the tangent vector to the equator, one has: 
$$
Q'(t)=(-\sqrt{a} \sin t, \sqrt{b} \cos t, 0), \ \   \langle Q'(t),Q'(t)\rangle=f^2(t).
$$ 
The geodesic corresponding to $(t,\alpha)$ has a tangent vector
$$
(-\sqrt{a}\sin t, \sqrt{b} \cos t, f(t) \tan\alpha).
$$
If the geodesic is time-like, we normalize the tangent vector so that its squared length is $-1$:
\begin{equation} \label{uvw}
(u,v,w)=\left( \frac{-\sqrt{a}\tau \sin t}{f(t)}, \frac{\sqrt{b}\tau\cos t}{f(t)}, \tau\tan\alpha \right).
\end{equation}

To obtain $J$, substitute $(u,v,w)$ from (\ref{uvw}) and 
\begin{equation} \label{xyz}
(x,y,z)=(\sqrt{a} \cos t, \sqrt{b} \sin t, 0)
\end{equation}
 to (\ref{Joa});  this yields the first formula (\ref{JandO}). 

Identify the cotangent and tangent bundles via the metric. Then the canonical symplectic form on the cotangent bundle becomes
$$
\omega= du\wedge dx+dv\wedge dy-dw\wedge dz.
$$
Formulas (\ref{uvw}) and (\ref{xyz}) describe a section of the tangent bundle, and the 
pull back of the symplectic structure $\omega$ is given by the second formula (\ref{JandO}). 
\proofend

As a consequence of Proposition \ref{onneg}, we define a natural  1-form on the space of null geodesics (this space consists of two disjoint circles corresponding to the right and left null geodesics). The set of null geodesics is given, in $(t,\alpha)$ coordinates, by $\alpha=\pi/4$ or $3\pi/4$, which corresponds to $\tau=\infty$. Note that both $\omega$ and $J$ blow up as one approaches the space of null geodesics, that is, as $\tau \to \infty$.  However, 
their ratio is well-defined.

\begin{lemma-def}\label{form}
The 1-form $h(t)dt$, given by the condition
that, in the limit $\tau\to\infty$,
\begin{equation} \label{divide}
d(J^{1/2})\wedge h(t)dt=\omega,
\end{equation}
is uniquelly defined. Explicitly, this equation holds for $h(t)$ 
given by
\begin{equation} \label{formh}
h(t)=const\cdot\frac{f(t)}{\sqrt{c+ f^2(t)}}=const\cdot\sqrt{\frac{a\sin^2t+b\cos^2t}{c+a\sin^2t+b\cos^2t}}
\end{equation}
with an appropriate constant factor.
\end{lemma-def}

\proof Up to a constant multiplier, one has:
$$
J^{1/2}=\tau \sqrt{c+f^2(t)} \left( 1+O\left( \frac{1}{\tau^2}\right)\right),
$$ 
hence 
$$
d J^{1/2} = \sqrt{c+f^2(t)}\,d\tau +g(t,\tau)\, dt+O\left( \frac{1}{\tau^2}\right)
$$
 for some function $g$. Therefore
$$
\left( d J^{1/2} \wedge \frac{f(t)}{\sqrt{c+ f^2(t)}}\, dt\right) =\omega +O\left( \frac{1}{\tau^2}\right),
$$
as needed. 
\proofend


\section{Poncelet-style closure theorem} \label{Poncelet}

\begin{figure}[hbtp]
\centering
\includegraphics[width=3.5in]{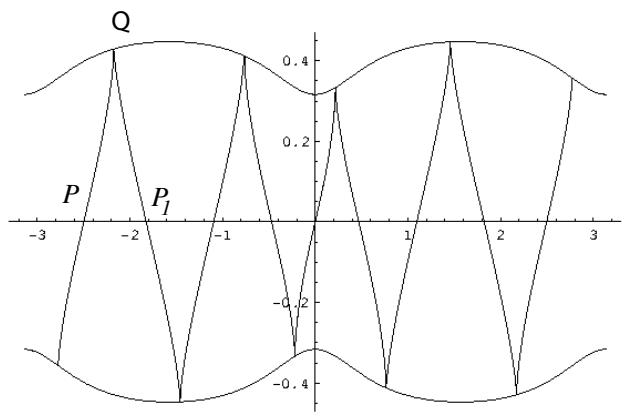}
\caption{The right null geodesic through $P$ and the left null geodesic 
through $P_1$ hit the Northern tropic at the same point $Q$.}
\label{Poncfig}
\end{figure}

Define a map $T$ of the equator to itself as follows. Given a point $P$ of the equator, consider the right null geodesic through $P$ until it intersects the Northern tropic 
at point $Q$, and then consider the left null geodesic through $Q$ until it intersects the equator at point $P_1$. Set:  $T(P)=P_1$, see figure \ref{Poncfig}.

\begin{theorem} \label{Poncthm}
Suppose that some point $P$ of the equator
is $k$-periodic, that is, $T^k(P)=P$ for some positive integer $k$.
Then  every  point of the equator is also $k$-periodic.
\end{theorem}

\proof
Consider the billiard inside the domain on the ellipsoid $M_0$  bounded by the curve $\Gamma_{\lambda}=M_0 \cap M_{\lambda}$. We claim that the respective billiard ball map  $F_{\lambda}$, acting on time-like geodesic segments, preserves the Joachimsthal integral $J$ considered as a function on the space of geodesics. A similar fact for Euclidean ellipsoids is well known, see, e.g., \cite{CS,Ve}.

Consider a time-like geodesic segment $\gamma_1$ on the ellipsoid, reflecting in $\Gamma_{\lambda}$ to another geodesic segment $\gamma_2$. According to Section \ref{intro}, the straight lines, tangent to $\gamma_1$ and $\gamma_2$, are tangent to  pseudo-confocal quadrics, say, $M_{\mu_1}$ and $M_{\mu_2}$. We want to show that $\mu_1=\mu_2$. Indeed, the billiard system in the ambient space with reflection in the quadric $M_{\lambda}$ is integrable and $M_{\lambda}$ is orthogonal to $M_0$. Therefore the incoming and the outgoing rays are tangent to the same pseudo-confocal quadric, see Section \ref{intro}, and thus $\mu_1=\mu_2$.

According to Section \ref{intro}, $F_{\lambda}$ preserves the symplectic structure $\omega$ on the space of time-like geodesics; it also preserves the integral $J$. Therefore the 1-form $h(t) dt$ on the space of null geodesics is also invariant under the action of the map $F_{\lambda}$ on the space of null geodesics. 

In the limit $\lambda \to 0$, the curve $\Gamma_{\lambda}$  becomes the tropic (\ref{tropic}), and the billiard ball map $F_{\lambda}$, restricted to null geodesics, gets identified with the map $T$. Hence $T$ preserves the 1-form $h(t) dt$. 

Finally, choose a cyclic coordinate $s$ on the equator so that $h(t) dt=ds$. In this coordinate, the map $T$ is a shift $s\mapsto s+c$. This map is $k$-periodic if and only if $kc\in\Z$.  This implies the statement of the theorem.
\proofend

\begin{problem} \label{Cayley}
{\rm It is interesting to find the relation on $a,b,c$, necessary and sufficient for 
the orbits of the map $T$ to close up after $k$ iterations and $r$ turns around the equator.  In the case of the Poncelet Porism, such conditions were found by Cayley, see \cite{G-H} for a modern treatment.
}
\end{problem}

\begin{remark} \label{otherP}
{\rm One also has Poncelet-style closure theorems for the geodesics in the polar caps and for space- or time-like geodesics in the equatorial belt reflecting from the tropics. Such results are similar to the ones known for Euclidean ellipsoid, see, e.g., \cite{CS,Ve}, and their proofs are similar to that of Theorem \ref{Poncthm} but simpler: the $T$-invariant 1-form $h(t) dt$ is obtained from a finite area form and a finite integral, not as a finite ratio of two infinite quantities.
}
\end{remark}


\section{Curvature of the ellipsoid and a geodesically equivalent Riemannian metric on it} \label{equi}

\subsection{Curvature of the ellipsoid}
The behavior of the geodesics in a polar cap  resembles that of the geodesics in the Poincar\'e disc model of the hyperbolic plane; this observation is explained by the following proposition.

\begin{proposition} \label{negcurv}
The Gauss curvature $K$ of the ellipsoid is negative; it is given by the formula:
$$
\frac{1}{K}=-abc\left( \frac{x^2}{a^2}+\frac{y^2}{b^2}-\frac{z^2}{c^2}\right)^2.
$$

\end{proposition}

\proof
Similarly to the Euclidean case, one can use the normal Gauss map to compute the Gauss curvature of a surface, see \cite{ON}. Consider the Northern polar cap. Normalize the normal vector as follows:
$$
N(x,y,z)=\frac{\left( \frac{x}{a},\frac{y}{b}, -\frac{z}{c} \right)}{\sqrt{\frac{z^2}{c^2}-\frac{x^2}{a^2}-\frac{y^2}{b^2}}},
$$
so that $\langle N,N\rangle =-1$. Thus the Gauss map sends the polar cap to
the upper sheet of the hyperboloid of two sheets, and this map reverses the orientation. Hence the Gauss curvature is negative.

Similarly, one considers the equatorial belt: the normal vector is then normalized to $\langle N,N\rangle =1$.

It is straightforward but tedious to compute the Gauss curvature, and we do not dwell on this. The computation can be simplified by the use of Mathematica.
\proofend

\subsection{Riemannian geodesically equivalent metric}

The next result follows from general constructions in \cite{M-T1,Ta3}; for completeness, we give a direct proof.  Two metrics are called {\it geodesically equivalent} if they have the same non-parameterized geodesics. 

\begin{proposition} \label{seconfm}
The metric on the ellipsoid (\ref{triell}), induced from the ambient Minkowski space, is geodesically equivalent to the Riemannian metric
\begin{equation} \label{other}
ds^2=\frac{\frac{dx^2}{a}+\frac{dy^2}{b}+\frac{dz^2}{c}}{\left| \frac{x^2}{a^2}+\frac{y^2}{b^2}-\frac{z^2}{c^2}\right|}.
\end{equation}
\end{proposition}

\proof
Let $P=(u,v,w)$ denote a tangent vector at point $Q=(x,y,z)$. Denote the diagonal matrix with the entries $(1/a,1/b,1/c)$ by $A$. To fix ideas, consider the equatorial belt.
The Lagrangian for the metric (\ref{other}) is
$$
L(P,Q)=\frac{\frac{u^2}{a}+\frac{v^2}{b}+\frac{w^2}{c}}{ \frac{x^2}{a^2}+\frac{y^2}{b^2}-\frac{z^2}{c^2}}=\frac{A(P)\cdot P}{f(Q)}
$$
with 
$$
f(Q)=\frac{x^2}{a^2}+\frac{y^2}{b^2}-\frac{z^2}{c^2}.
$$
The Euler-Lagrange equation for a geodesic $Q(t)$  with $\dot Q =P$ is
$$
L_{PP} (\ddot Q) + L_{PQ} (P) - L_Q = \lambda A(Q)
$$
where $A(Q)$ is a Euclidean normal to the ellipsoid at point $Q$ and $\lambda$ is a Lagrange multiplier; here $L_{PP}$ and $L_{PQ}$ are the matrices of the second partial derivatives and $L_Q$ is the gradient vector.

One easily computes:
$$
L_{PP}=\frac{2}{f(Q)} A,\ L_{PQ}=-\frac{2}{f^2(Q)} A(P)\otimes \nabla f(Q),\ 
L_Q=-\frac{A(P)\cdot P}{f^2(Q)} \nabla f(Q),
$$
and the Euler-Lagrange equation is rewritten as
\begin{equation} \label{EL}
A(\ddot Q)-\frac{P\cdot \nabla f(Q)}{f(Q)} A(P) + \frac{A(P)\cdot P}{2f(Q)} \nabla f(Q) =\lambda A(Q).
\end{equation}
To find the Lagrange multiplier, dot multiply equation (\ref{EL}) by $Q$. One has: 
$$
A(Q)\cdot Q=1,\ A(P)\cdot Q=P\cdot A(Q)=0,\ \nabla f(Q)\cdot Q=2f(Q),
$$
and hence 
$$
\lambda= A(\ddot Q)\cdot Q + A(P)\cdot P = \frac{d({A(P)\cdot Q})}{dt}=0.
$$
Thus (\ref{EL}) implies that the acceleration $\ddot Q$ lies in the plane spanned by the velocity vector $P$ and the vector $A^{-1} (\nabla f(Q))=2N(Q)$, the Minkowski  normal to the ellipsoid at point $Q$. It follows that $Q(t)$ is a reparameterized geodesic of the restriction of the ambient metric on the ellipsoid.

The case of the polar caps is similar.
\proofend  

\begin{remark} \label{ratio}
{\rm Note that the metric (\ref{other}) is given by the ratio of the two factors whose product is the Joachimsthal  integral (\ref{Joa}).
}
\end{remark}


\section{Rigidity of ellipses as integrable Lorentz billiard curves} \label{rig}

In Section \ref{pancakes}, we considered the limit of the ellipsoid (\ref{triell}) as $b\to 0$. In this limit, the ellipsoid becomes an ellipse in the Lorentz plane. This ellipse is foliated by two parallel families of null lines, and the map $T$, defined at the beginning of Section \ref{Poncelet}, takes oriented lines in one null direction to oriented null lines in the other one.

More generally, consider an oval (closed smooth strictly convex curve) $\gamma$ in the Lorentz plane and the billiard system inside it. The billiard reflection takes one null direction to the other. Thus, restricted to the null lines, the billiard ball map is identified with the  circle map $T:\gamma \to \gamma$ described as follows. The intersections with the lines in the first direction define an involution on $\gamma$, and likewise for the second direction. The composition of these involutions is the map $T$. 

\begin{figure}[hbtp]
\centering
\includegraphics[width=2in]{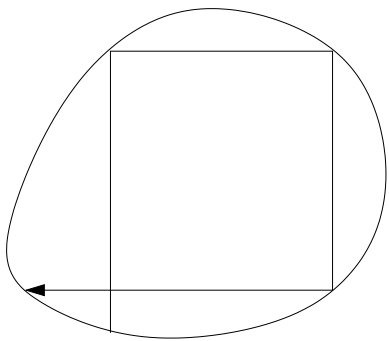}
\caption{A map of an oval}
\label{arn}
\end{figure}

The circle map $T$ is depicted in figure \ref{arn}; it was previously considered in different contexts: in relation to Hilbert's 13th problem \cite{Fest}; the Sobolev equation, approximately describing fluid oscillations in a fast rotating tank, \cite{Sob1}; and the theory of Lorentz surfaces, where the rotation number of $T$ provides a continuous invariant  of the conformal class of a Lorentz disc \cite{Wei}.

The map $T:\gamma \to \gamma$ depends only on the choice of two (null) directions, say, $u$ and $v$; thus we denote it by $T_{(u,v)}$. If $\gamma$ is an ellipse then $T_{(u,v)}$ is conjugated to a rotation. This is a limit case of Theorem \ref{Poncthm} but it is also easily proved directly:  an affine transformation  takes $\gamma$ to a circle, and then $T_{(u,v)}$ becomes the rotation through the angle twice that between the two respective  directions.

\begin{conjecture}
Let $\gamma$ be a plane oval such that, for every pair of directions $u$ and $v$, the map $T_{(u,v)}$ is conjugated to a circle rotation. Then $\gamma$ is an ellipse.
\end{conjecture}

We prove a weaker result. 

\begin{theorem} \label{char}
Let $\gamma(t)$ be a parameterized plane oval with the property that, for every pair of directions $u$ and $v$, the map $T$ is a translation in the variable $t$, i.e., there is a constant $c(u,v)$ such that  $T_{(u,v)} (\gamma(t))=\gamma(t+c(u,v))$. Then $\gamma$ is an ellipse.
\end{theorem}

\proof
Without loss of generality, assume that $\gamma$ is parameterized by $t\in\R/2\pi\Z$.

Choose a direction $u$. There is a unique maximal chord of $\gamma$  in direction $u$, say, $AB$; such a chord is called an {\it affine diameter}.  Since $AB$ is maximal, the tangent lines to $\gamma$ at $A$ and $B$ are parallel; let $v$ be their direction. It is convenient to apply an affine transformation that makes $u$ horizontal and $v$ vertical.

 Let $C$ and $D$  be the points at which the tangent line to $\gamma$ are horizontal. We claim that the chord $CD$ is vertical.

\begin{figure}[hbtp]
\centering
\includegraphics[width=2in]{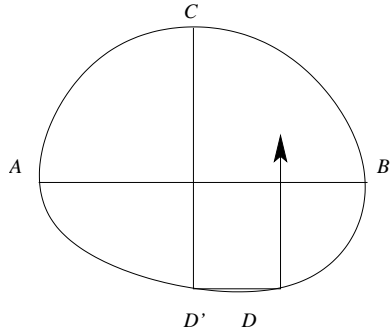}
\caption{Proving that $D=D'$}
\label{skew}
\end{figure}

Indeed, consider the map $T_{(u,v)}$. Then $T^2(A)=A$. Since $T_{(u,v)}$ is conjugated to  a rotation, every point of $\gamma$ is 2-periodic.  Let $CD'$ be the vertical chord of $\gamma$. Since $T_{(u,v)}^2(C)=C$, the  direction of $\gamma$ at $D'$ is horizontal, and hence $D'=D$, see figure \ref{skew}.

Thus, for every direction $u$, there exists a conjugate direction $v$: there are affine diameters $AB$ and $CD$ of the curve $\gamma$, having directions $u$ and $v$, such that the tangent lines at $A$ and $B$ are parallel to $CD$, and the tangent lines at $C$ and $D$ are parallel to $AB$ (this makes $\gamma$ a Blaschke $P$-curve \cite{Bla}). 

Furthermore, if $(u,v)$ is a pair of conjugate directions then every point of $\gamma$ is a vertex of an inscribed parallelogram with the sides parallel to $u$ and $v$, and whose diagonals correspond to a 2-periodic orbit of $T_{(u,v)}$.  The  opposite vertices of these parallelograms are 2-periodic points $\gamma(t)$ and $\gamma(t+\pi)$.

\begin{figure}[hbtp]
\centering
\includegraphics[width=1.2in]{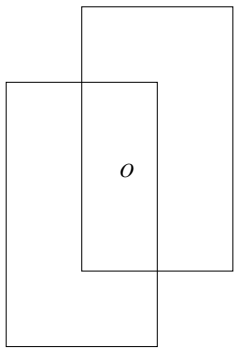}
\caption{Inscribed rectangles}
\label{rect}
\end{figure}

Assume that $\gamma$ is centrally symmetric with respect to the origin $O$. Then each rectangle is centered at $O$ -- otherwise $\gamma$ is not convex, see figure \ref{rect}.
 It follows that the mid-points of the horizontal sides of the rectangles lie on the affine diameter $CD$, and the mid-point of the vertical sides lie on $AB$. This is true for every pair of conjugate directions, hence $\gamma$ admits affine line symmetry for every direction. This is a characteristic property of ellipses \cite{Ber}, and thus $\gamma$ is an ellipse.

\begin{figure}[hbtp]
\centering
\includegraphics[width=1.8in]{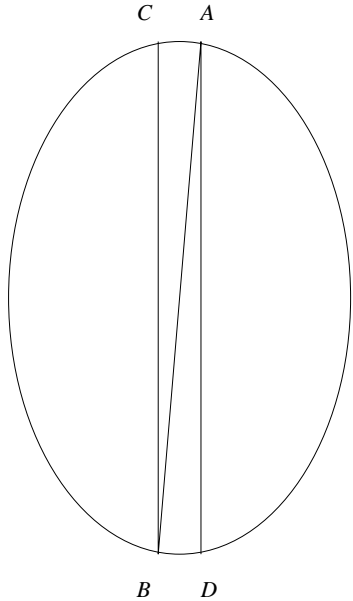}
\caption{Proving central symmetry of $\gamma$}
\label{sym}
\end{figure}

It remains to prove that $\gamma$ is centrally symmetric. Consider two directions, $u$ and $v$, making an infinitesimal angle $\varepsilon$. Consider figure \ref{sym} in which $AB$ is an affine diameter in the direction $u$, and the lines $BC$ and $AD$ have the direction $v$. One has: $T_{(u,v)} (A)=C, T_{(u,v)} (B)= D$. Let $A=\gamma(t),  C=\gamma(t+\delta)$. Then $B=\gamma(t+\pi), D=\gamma(t+\pi+\delta)$.
Since  the opposite sides of  the infinitesimal quadrilateral $ACBD$ are parallel, it is 
 a parallelogram and $AC=BD$. It follows that $\gamma'(t)=-\gamma'(t+\pi)$, and hence $\gamma(t)+\gamma(t+\pi)$ is a constant vector. Thus $\gamma$ is centrally symmetric.
\proofend

\bigskip

{\bf Acknowledgments}. We are grateful to V. Chernov, D. Fuchs and P. Lee for stimulating discussions. The second and the third authors were partially supported by NSERC and NSF grants, respectively.

\end{document}